\documentclass[english,12pt]{article}

\usepackage{amsmath,amsbsy,amscd,amssymb,graphicx,color}

\def \ts{\textstyle}

\def \varep{\varepsilon}

\def \ul{\underline}

\def \gs{\geqslant}
\def \ls{\leqslant}
\def \part{\partial}

\begin{document}

\title{Controlling Chaotic Maps by \\ Feedback Control Modulation} 

\author{Roberta Hansen \& Graciela Adriana Gonz\'alez \\ \\ \it \small Departamento de Matem\'atica, Facultad de Ingenier\'{\i}a, Universidad de Buenos Aires, Argentina, \\ \it \small Consejo Nacional de Investigaciones Cient\'{\i}ficas y T\'ecnicas, Argentina. \\ \tt \small rhansen@fi.uba.ar, ggonzal@fi.uba.ar
}
\date{}
\maketitle

\begin{abstract}
This paper deals whith the stabilization of any UPO of a chaotic map by modulation of a control parameter. It concentrates on proportional and delayed feedback control methods. Alternative types of these methods are proposed and their achievments are investigated analyticaly and numerically.
\end{abstract}

{\bf keywords:} chaos control, logistic map, proportional feedback control, delayed feedback control, parameter modulation. 

37N35 - 93D15

{\thispagestyle{empty}} 

\section{Introduction}

Chaotic behavior is a very interesting nonlinear phenomena, but in many practical situations it is desirable to be avoided, for example, because it restricts the operating range of electronic or mechanics devices. Moreover, this goal must be achieved with the only help of tiny perturbations properly chosen \cite{Rom, Shinb, Bocca}.

The seminal idea of Ott, Grebogy and Yorke \cite{OGY} turned the presence of chaos into an advantage. Indeed, the system may be stabilized in a particular unstable periodic orbit (UPO) embedded within a strange attractor. When the trajectory is close to a desired UPO, a small time-dependent feedback perturbation is applied to some accessible parameter or variable system. The periodic orbit is preserved, but its stability is modified keeping the trajectory to stay close to the UPO. This control strategy is known as the OGY method.

A simple proportional feedback (SPF) control method, basically consists in a perturbation proportional to the difference between the current system value and an unstable fixed point (or an UPO's component). The OGY method (\cite{OGY}) is a particular case of a SPF. It is well known that it can control chaos for some 1-dimensional maps \cite{G, Bocca}. The usefulness of SPF control algorithms for 1-dimensional maps arises when controlling chaos in highly dissipative systems; as there are cases in which it loses validity because Poincar\'e sections changes in each iteration, a recursive proportional feedback (RPF) has been proposed (see \cite{WCHO} and references therein). Both, SPF and RPF methods require the exact knowledge of the UPO to be stabilized, and the linearized dynamics about it. This is not always at one's disposal in real-world implementations. Moreover, the OGY method results very sensitive to nonlinearities and fluctuations of external noise, mainly for large orbit periods. With the aim to overcome these limitations, Pyragas \cite{P} proposed an alternative method, based on a self-controlling delayed feedback, which, as the OGY, does not modify the original UPO, but it does not depend explicitely on it. Its discrete-time version results in a perturbation which is proportional to the difference between the current system value and a previous one. This delayed feedback control (DFC) successfully controls chaotic behavior in a variety of experiments \cite{HT} (and its references). However, the stabilization capability of the DFC may be weaker than the OGY's. The Pyragas method is renamed as ``time-delay-autosynchronization'' (TDAS) in \cite{Socolar}, where an extended version of it is proposed which uses information from many previous states, and so, improving on stabilization objectives achievement; this extension is refered as ETDAS in \cite{Socolar} and EDFC in \cite{Ko}. It is worth to mention that there exist other proposals in the literature oriented to overcome these drawbacks, among which, the oscillating approach \cite{Morgul}, the predictive one \cite{UY}, and the adaptive version \cite{PChagas} are, perhaps, the most interesting. 

Assume that a discrete-time dynamical system is given by
\begin{equation} 
\label{discrete-map}
x_{k+1}=f(x_k,r)\, ,
\end{equation}
where $r$ is a scalar parameter, and that for $r\!=\!r_0$ the system developes chaotic behavior, having an infinite number of UPOs embedded within a strange attractor, $\cal{A}$, which is included in its basin of attraction, $\cal{B}$. The dynamics in the chaotic attractor is ergodic, meaning that for almost all initial condition in $\cal{B}$, a system trajectory visits any small neighborhood of every point of the UPO's.

The logistic map is the prototypical discrete-time dynamical system with simplest algebraic equation and exhibiting all ingredients of chaos, and so extensively studied: 
\begin{equation} 
\label{logistica}
x_{k+1}=f(x_k,r)=r\,x_k\,(1\!-\!x_k)\, ,
\end{equation}
where $x\!\in\![0,1]$ and $r$ is a control parameter. This map develops a cascade of period--doubling bifurcation that accumulates at $r\!=\!r_{\infty}\!\approx\!3.569$, where the chaotic regime starts with the presence of chaotic attractors, until $r\!=\!4$, from which the motion becomes unbounded. For $r_{\infty}\!<\!r_0\!\ls\!4$ there is one chaotic attractor with ${\cal{B}}\!=\![0,1]$

A widely used technique for controlling chaotic behavior is modulation: by judiciously varying control parameters, system trajectory is driven into a desirable dynamical state. We want trajectories resulting from a randomly chosen initial condition $x_0$ to be as close as possible to a $m$-period UPO, $\{p_1,p_2,\hdots,p_m\}$. Suppose, that the parameter $r$ can be finely tuned in a small range around $r_0$, namely, we allow $r$ to vary in the range $(r_0\!-\!\delta,r_0\!+\!\delta)$, where $0\!<\!\delta\!\ll\!1$. Then, our objective is to stabilize Eq. (\ref{discrete-map}) at this UPO by feedback control modulation. Namely, the (control) parameter is affected by a control $u_k\!=\!u(x_k)$, so
\begin{equation}\label{map-control}
x_{k+1}\!=\!f(x_k,r+u_k), 
\end{equation}
under the requirement that the adjusted parameter remains within a range for which the system is chaotic in the absence of perturbations.

Published works on controlling chaos by modulation are mostly based on experimental or numerical arguments, while analytical approaches use to concentrate on a fixed point stabilization problem \cite{Socolar2}. Many results are developed on the logistic map (\cite{ES, NDBR} and references therein). DFC and SPF controls applied to 1-dimensional maps are confronted in \cite{HT} through linear stability analysis and by implementation in an analog electric circuit. 

In this work, we investigate analytically and numerically the approach to stabilize any UPO of the chaotic map (\ref{discrete-map}) by modulation of the parameter $r$. We explore the possibility of relaxing  some requirements of the well known OGY and Pyragas methods, and we also propose alternative types of SPF and DFC methods. The logistic map will be repeatedly used to illustrate the implementation of the methods and so appreciate their control performance. The parameter modulation is stated through the following controlled system
\begin{equation}\label{param_modul}
x_{k+1}=[r_0+u_k]\, x_k\,(1\!-\!x_k)\ .
\end{equation}
Throughout the work, the problem statement for the general case of Eq. (\ref{discrete-map}) and its theoretical approaches will be dealt with. By simplicity, first it will be developed for 1-dimensional maps and finally, its extension to the $n$-dimensional case will be considered.







\section{Controlling Chaos by SPF Modulation}

\subsection{The OGY method}
 
We briefly give the ideas of this method outlined in \cite{OGY}. Assume that at time $k$, the trajectory falls into the neighborhood of the component $p_i$. The linearized dynamics of Eq. (\ref{discrete-map}) about $p_i$ and $r_0$ is:
\begin{equation}\label{linearlogistic}
x_{k+1}\!-\!p_{i+1}\ =f_x(p_i,r_0)\,(x_k\!-\!p_i)\,+\,f_r(p_i,r_0)\,(r_k\!-\!r_0)\,.
\end{equation}
To force the system towards the UPO, we set $x_{k+1}\!-\!p_{i+1}\!=\!0$, and so we have, from (\ref{linearlogistic}): 
\begin{equation}
\label{delta-r}
u^i_k=(\Delta r)_k=r_k\!-\!r_0=-\frac{f_x(p_i,r_0)}{f_r(p_i,r_0)}\,(x_k\!-\!p_i)= \alpha_i\,(x_k\!-\!p_i)\, .
\end{equation}
Eq. (\ref{delta-r}) holds only when $|x_k\!-\!p_i|\!\ls\!\varep\!\ll\!\!1$, hence the required parameter perturbation $(\Delta r)_k$ is small, and the maximum parameter perturbation $\delta$ is proportional to $\varep$ (with factor $\alpha_i$). When the trajectory is outside the $\varep$-neighborhood of $p_i$, the perturbation is not applied, and the system evolves at its nominal chaotic parameter $r_0$. Then, for given $\varep,\delta\!>\!0$, the control algorithm is
\begin{equation}\label{controlOGY}
u^i_k=\left\{\!\begin{array}{cc}
\alpha_i\,(x_k\!-\!p_i) & \ {\rm if }\quad |x_k\!-\!p_i|\!\ls\!\varep \, ,\\
  0   & \quad {\rm otherwise}\, .
\end{array}\right.
\end{equation}
This dynamics preserves the $m$-UPO by a parameter perturbation of SPF type, i.e., $u^i_k\!\propto\!|x_k\!-\!p_i|$. 

In the case of the Eq. (\ref{logistica}), $f_x(p_i,r_0)\!=\!(1\!-\!2p_i)$ and $f_r(p_i,r_0)\!=\!p_i(1\!-\!p_i)$, so $\alpha_i\!=\!\frac{\ts r_0(2\,p_i\!-\!1)}{\ts p_i(1\!-\!p_i)}$. As the modulation of the parameter $r$ should keep the dynamics to remain globally bounded, it must be $r_0\!+\!|u_k|\!\ls\!4$, $r_0\!\in\!(3.569,4)$, and hence, for $\varep\!>\!0$, $|\alpha_i\,\varep|\!\ls\!\delta\!\ls\!4\!-\!r_0$, or $\varep\!\ls\!\delta/|\alpha_i|$, for all $1\!\ls\!i\!\ls\!m$.
\subsection{A slight modification of the OGY method}

The control in Eq. (\ref{controlOGY}) also depends on which $p_i$ was selected to stabilize the dynamics on it, namely it is turned on only at the end of each oscillation. This fact makes this control very sensitive to nonlinearities and to fluctuation of external noises that are common in real implementations, mainly for large values of the orbit's period, $m$. Besides, for every $i$, different control gain and waiting time come out, yielding to different performances depending on each case.

\begin{figure}[h]
\hspace{-3cm}
\includegraphics[width=21cm,height=15cm]{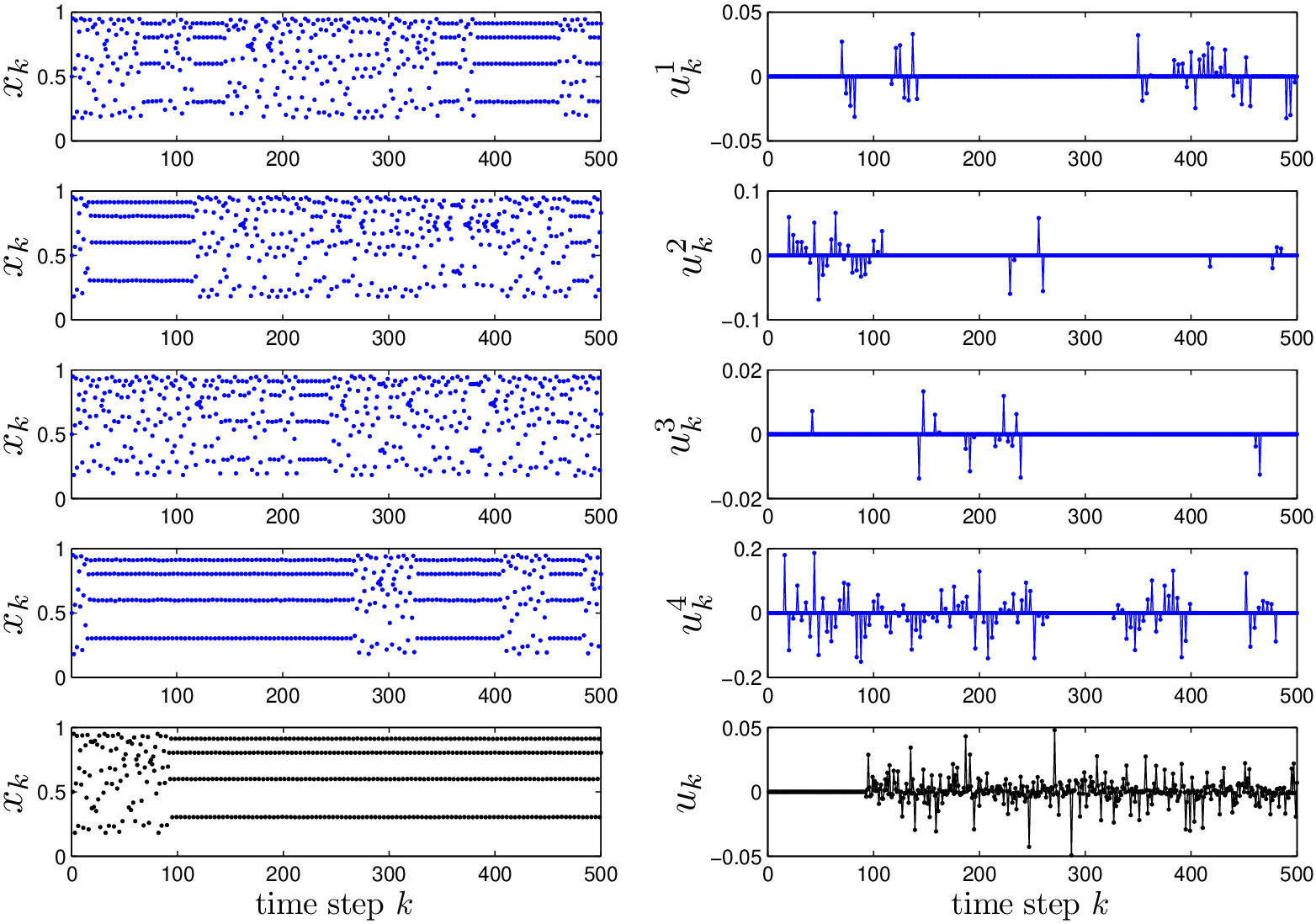}
\caption{\small (blue) States $x_k$, and controls $u^i_k,i\!=\!1$--$4$ of Eq. (\ref{controlOGY}) with $\varep\!=\!0.005$, applied to stabilize $4$-UPO of Eq. (\ref{logistica}) for $r\!=\!3.8$, $\{p_1\!\approx\!0.3,p_2\!\approx\!0.8,p_3\!\approx\!0.6,p_4\!\approx\!0.91\}$, under the effect of additive noise modeled by $5\!\times\!10^{-4}\sigma_k$, $\sigma_k\!\sim\!N(0,1)$, i.c. $x_0\!=\!0.5,\varep\!=\!0.005$. (black) The same, but applying the control (\ref{switching}).}
\label{cont-alf-ruido}
\end{figure}

In order to improve these features, we first propose a simple ``switching control'', that works as the OGY, but forcing the trajectory to keep close to the $m$--period UPO, by applying the perturbation $u_k^i$ for all $1\!\ls\!i\!\ls\!m$, i.e., each time the trajectory is close to any component $p_i$. As a first and natural consequence, a net reduction on waiting time will be obtained. The control algorithm is
\begin{equation}\label{switching}
u_k=\left\{\!\begin{array}{cc}
\alpha_i\,(x_k\!-\!p_i) & \quad {\rm if }\quad |x_k\!-\!p_i|\!\ls\!\varep,\ 1\!\ls\!i\!\ls\!m \, ,\\
  0   & \quad {\rm otherwise}\, ,
\end{array}\right.
\end{equation}
provided that $\varep<\frac{\ts |p_i\!-\!p_j|}{\ts 2},\ \forall\, i\!\neq\! j$. As a second consequence, this strategy of control displays a notably better performance than control (\ref{controlOGY}) in presence of external noise. This is numerically verified in the logistic for which the the response of applying the control (\ref{switching}) about a 4-period UPO of map (\ref{logistica}) with $r_0\!=\!3.8$, and the controls $u^i_k$ of (\ref{controlOGY}) for each $i\!=\!1,2,3,4$, in presence of noise, are compared in Fig. \ref{cont-alf-ruido}.

\subsection{Improving SPF modulation}

Here we propose to perturb the parameter $r_0$ with a SPF modulation selected from a set of control laws, similar to (\ref{switching}), but replacing the fixed gain $\alpha_i$ by coefficients $\beta_i$ adequately chosen, 
\begin{equation}\label{controlbetas}
u_k=\left\{\!\begin{array}{cc}
\beta_i\,(x_k\!-\!p_i) & \quad {\rm if }\quad |x_k\!-\!p_i|\!\ls\!\varep,\  1\!\ls\!i\!\ls\!m \, ,\\
  0   & \quad {\rm otherwise}\, ,
\end{array}\right.
\end{equation}
provided that $\varep<\frac{\ts |p_i\!-\!p_j|}{\ts 2},\ \forall\, i\!\neq\! j$. This dynamics also preserves the $m$-UPO. The linear stability criterion for $p_i$ to be an asymptotically stable (a.s.) point, states the condition for $\beta_i$. Introducing $f_{\beta_i}(x)=f(x,r+\beta_i\,(x-p_i))$, let us assume that there exists a range $(\beta_i^{\inf},\beta_i^{\sup})$ of $\beta_i$ values such that $\left|f_{\beta_i}'(p_i)\right|\!<\!1,\ 1\!\ls\! i\!\ls\! m$. 
For the logistic map (\ref{logistica}), $f_{\beta_i}(x)\!=\!r+\beta_i\,(x-p_i)\,x\,(1\!-\!x)$, being
$$\beta_i^{\inf}= -\, \frac{1+r_0(1\!-\!2p_i)}{p_i(1\!-\!p_i)}\qquad {\rm and}\qquad  \beta_i^{\sup}=\frac{1\!-\! r_0(1\!-\!2p_i)}{p_i(1\!-\!p_i)}$$
For example, for $r_0\!=\!3.8$ the $\beta_i$ must verify $|\beta_i\,\varep|\!\ls\!\delta\!\ls\!0.2$, $\forall\, i$, to ensure the desired bound on the control effort and a globally bounded dynamics. Under the conditions stated on the control coefficients, convergence of the algorithm is formally proven.
\medskip

{\bf Proposition 1:} {\it Let the controlled system (\ref{map-control})--(\ref{controlbetas}) with $\beta_i\in (\beta_i^{\inf},\beta_i^{\sup})$, $1\!\ls\! i\!\ls\! m$ and $\beta\!=\!\max\limits_{1\ls i\ls m}\{|\beta_i|\}$. Then, there exists $\varep_0$, $0\!<\!\varep_0\!<\!\min\limits_{i\neq j}\Big\{\frac{\ts |p_i\!-\!p_j|}{\ts 2}; \frac{\ts \delta}{\ts \beta}\Big\}$ such that for all $\varep$, $0\!<\!\varep\!<\!\varep_0$, and for almost every initial condition $x_0\!\in\!\cal{B}$, it verifies $|u_k|\ls \delta \  \forall\, k$, and $(x_k)_{k\gs 1}$ converges to the $m$-periodic orbit $\{p_1,\hdots,p_m\}$.}
\medskip

\ul{Proof:} 

Let $\varep'_0=\min\limits_{i\neq\j}\Big\{\frac{\ts |p_i\!-\!p_j|}{\ts 2}\Big\}$. For all $i$, there exists $\sigma>0$, such that $\left|f_{\beta_i}'(p_i)\right|<\sigma<1$, and, for every $i$, there exists $0<\varep_i<\varep_0'$ such that $\left|f_{\beta_i}'(x)\right|\ls \sigma$, for $x\in (p_i-\varep_i,p_i+\varep_i)$. Let $\varep_0\!=\!\min\limits_{1\ls i\ls m}\varep_i$. 

Fix $\varep\!<\!\varep_0$, for almost every $x_0\!\in\!\cal{B}$ the ergodicity of the uncontrolled system guarantees the existence of $k_0=k(x_0,\varep)$, such that $|x_{k_0}-p_i|\ls\varep$ for some $i$, then 
$$|x_{k_0+1}-p_{i+1}|= |f_{\beta_i}(x_{k_0})-f_{\beta_i}(p_i)|= \left|f_{\beta_i}'(\xi_i)\right||x_{k_0}-p_i|\ls\sigma\varep<\varep$$
for $\xi_i\in (p_i-\varep_i,p_i+\varep_i)$, ($p_i\equiv p_{i\,{\rm mod}(m)}$). Applying successively the algorithm of Eq. (\ref{controlbetas}), we obtain
$$|x_{k_0+n}-p_{i+n}|\ls\sigma^n\varep<\varep,\quad \forall\, n\gs 1,$$
and the thesis follows. \hfill{$\square$}
\medskip

The existence of a range $(\beta_i^{\inf},\beta_i^{\sup})$ states the robustness of method of the Eq. (\ref{controlbetas}), and in particular of the Eq. (\ref{switching}) (note that for $\alpha_i\!=\!\frac{\ts \beta_i^{\inf}\!+\!\beta_i^{\sup}}{\ts 2}$, $\ f'_{\alpha_i}(p_i)\!=\!0\ $). As $\beta^{\inf}\!<\!\alpha$, by means of (\ref{controlbetas}) the objective is fulfilled with smaller values for the control gain, or else, for the same control effort, $\delta$, a greater $\varep$ is allowed, improving the waiting time to active the control.
As an example, we take the logistic map and its unstable fixed point $p\!=\!1\!-\!\frac{1}{r_0}$, requiring $\delta\!=\!0.2$ as bound on control effort. In Figure 2(a), it is appreciated both: the neat reduction of control effort by changing the control gain $\alpha$ by $\beta$, and an increase in the convergence time once the control is turned on. Figure 2(b) shows, for the similar control effort, the reduction in the waiting time explained above.
\medskip

Two remarks may be pointed out on (\ref{controlbetas}):
\begin{enumerate} \itemsep -0.2cm
\item[i)] verification of $\varep-$nearness becomes superfluous once $k_0$ is detected,
\item[ii)] the choosing of $\beta_i$ implies (but it is not equivalent to!)
\begin{equation}\label{prod_deriv}
\big|f'_{\beta_1}(p_1)f'_{\beta_2}(p_2)\cdots f'_{\beta_m}(p_m)\big|<1 \ .
\end{equation}
\end{enumerate}
These observations yield to a generalization of (\ref{controlbetas}) which, besides, includes (\ref{controlOGY}) as a particular case. This control is of the form
\begin{equation}\label{control_mod}
u_k=\left\{\!\begin{array}{cc}
0 &  \quad k<k_0\, ,\\
\beta_{(i+k-k_0)\,{\rm mod}(m)}(x_k-p_{(i+k-k_0)\,{\rm mod}(m)})   & \quad k\gs k_0\, ,
\end{array}\right.
\end{equation} 
where $\varep\!<\!\min\limits_{i\neq j}\Big\{\frac{\ts |p_i\!-\!p_j|}{\ts 2}\Big\},\ k_0\!=\!k_0(x_0,\varep)\!=\!\min\{k\!\gs\!0: |x_k\!-\!p_i|\ls \varep,\ {\rm for\ some}\ i,\ 1\!\ls\! i\!\ls\! m\}$. Note that $k_0$ exists by ergodicity.

\medskip

{\bf Proposition 2:} {\it Let the controlled system (\ref{map-control})-(\ref{control_mod}) with $\beta_j,\ 1\!\ls\! j\!\ls\! m$, for which (\ref{prod_deriv}) is valid, and $\beta\!=\!\max\limits_{1\ls j\ls m}|\beta_j|$. Then, there exists $\varep_0$, $0\!<\!\varep_0\!<\!\min\limits_{i\neq j}\Big\{\frac{\ts |p_i\!-\!p_j|}{\ts 2};\frac{\ts\delta}{\ts \beta}\Big\}$ such that for all $\varep$, $0\!<\!\varep\!<\!\varep_0$, and for almost every initial condition $x_0\!\in\!\cal{B}$, it verifies $|u_k|\!<\!\delta$, and $(x_k)_{k\gs 1}$ converges to the $m$-periodic orbit $\{p_1,\hdots,p_m\}$.} 
\medskip

\ul{Proof:} Let $\varep'_0\!=\!\min\limits_{i\neq j}\Big\{\frac{\ts |p_i\!-\!p_j|}{\ts 2}\Big\}$, and $0\!<\sigma\!<\!1$ such that
\begin{equation}\label{prod_deriv_sigma}
\big|f'_{\beta_1}(p_1)f'_{\beta_2}(p_2)\cdots f'_{\beta_m}(p_m)\big|<\sigma<1
\end{equation}
As (\ref{prod_deriv_sigma}) may be re-written
$$\big|f'_{\beta_i}(p_1)f'_{\beta_2}\big(f_{\beta_1}(p_1)\big)\cdots f'_{\beta_m}\big(f_{\beta_{m-1}}\big(\cdots\big(f_{\beta_1}(p_1)\big)\big)\big)\big|<\sigma$$
there exists $\varep_1,\ 0\!<\!\varep_1\!<\!\varep_0'$ such that
$$\big|f'_{\beta_i}(x)f'_{\beta_2}\big(f_{\beta_1}(x)\big)\cdots f'_{\beta_m}\big(f_{\beta_{m-1}}\big(\cdots\big(f_{\beta_1}(x)\big)\big)\big)\big|\ls\sigma$$
for $x\!\in\! (p_1\!-\!\varep_1,p_1+\varep_1)$. Proceeding analogously  for $p_2,\hdots,p_m$, the values $\varep_2,\hdots,\varep_m$ arise. Let $\varep_0\!=\!\min\limits_{1\ls j\ls m}\varep_j$, and fix $\varep\!<\!\varep_0$. Without loss of generality, we assume that the condition $|x_k\!-\!p_i|\!\ls\!\varep$ is verified, for the first time, by $i\!=\!1$. Then,
\begin{eqnarray*}
|x_{k_0+m}-p_1|&=&\big| f_{\beta_m}\circ f_{\beta_{m-1}}\circ\cdots\circ f_{\beta_1}(x_k)- f_{\beta_m}\circ f_{\beta_{m-1}}\circ\cdots\circ f_{\beta_1}(p_1)\big|\\
               &=& \big| \big( f_{\beta_m}\circ f_{\beta_{m-1}}\circ\cdots\circ f_{\beta_1}\big)'(\xi_1)\big| \, |x_k-p_1| \\
               &\ls& \sigma\, \varep < \varep
\end{eqnarray*}
for $\xi_1\!\in\!(p_1-\varep_1,p_1+\varep_1)$, and so following
\begin{equation}\label{xxx}
|x_{k_0+\ell m}-p_1|\ls \sigma^{\ell}\varep < \varep,\ \ell\!\gs\!1.
\end{equation}
By the same arguments, it is obtained that
\begin{equation*}
|x_{k_0+j+m}-p_{j+1}|\ls \sigma\,\varep < \varep\quad {\rm and}\quad |x_{k_0+j+\ell m}-p_{j+1}|\ls \sigma^{\ell}\varep < \varep,\ j\!=\!1,\hdots,m\!-\!1.
\end{equation*}
Hence, control bounds and convergence to the periodic orbit follow. \hfill{$\square$}

\section{Controlling Chaos by DFC Modulation}
\subsection{The Pyragas method}

In \cite{P} Pyragas proposed to stabilize the logistic map (\ref{logistica}) to its unstable fixed point, $p\!=\!1\!-\!\frac{\ts 1}{\ts r_0}$, by additive forcing in the form of one-time delay linear perturbation $\gamma(x_k\!-\!x_{k-1})$. This perturbation can be useful in applications were the fixed point is unknown or may drift. Here we use the same DFC but modulating the parameter $r$. Taking care of control bounds we propose
\begin{equation}\label{DFCptofijo}
u_k=\left\{\!\begin{array}{cc}
\gamma\,(x_k\!-\!x_{k-1}) & \ {\rm if }\quad \sqrt{|x_k\!-\!p|^2\!+\!|x_{k-1}\!-\!p|^2} \ls \frac{\ts \varep}{\sqrt{2}} \, ,\\
  0   & \quad {\rm otherwise}\, ,
\end{array}\right.
\end{equation}
Control (\ref{DFCptofijo}) vanishes when the system (\ref{map-control}) state attains $p$, so the fixed point is preserved. Let us see how it works in the logistic case. The controlled system (\ref{param_modul})-(\ref{DFCptofijo}) yields to the two dimensional dynamical system: 
\begin{equation*}
\left\{\begin{array}{l}
x_{k+1}=[r_0\!+\!\gamma\,(x_k\!-\!y_k)]x_k\,(1\!-\!x_k),\\
y_{k+1}=x_k,
\end{array}\right.
\end{equation*}
which has $P\!=\!(p,p)$ as a fixed point. The Jacobian matrix at $P$, $J(P)$, is a $(2\!\times\!2)$-{\em companion matrix} of the form
\begin{equation*}\left(
\begin{array}{cc} a_{11} & a_{12}\\ 1 & 0
\end{array}\right)
\end{equation*}
where $a_{11}=\frac{\gamma}{r_0}\big(1-\frac{1}{r_0}\big)+2-r_0$, and $a_{12}=-\frac{\gamma}{r_0}\big(1-\frac{1}{r_0}\big)$. The necessary and sufficient conditions for $P$ to be a.s. \cite{K}, yields to 
$$\gamma^{\inf}=\frac{r_0^2\,(r_0\!-\!3)}{2\,(r_0\!-\!1)}<\gamma<\gamma^{\sup}=\frac{r_0^2}{r_0\!-\!1} \ .$$
Fixing a $\gamma$ in this range, and choosing an adequate $\varep$ to assure the control effort to be bounded by $\delta$, the convergence of the trajectories to $p$ is obtained.

The resulting control performance by applying (\ref{DFCptofijo}), may be comparable, or even better, than the one obtained by applying (\ref{controlbetas}) if adequate coefficients are chosen (see Figures 3(a) and (b)).

\subsection{Improving DFC modulation} 
Our proposal of an extension of (\ref{DFCptofijo}) to a $m$-UPO, $\big\{p_1,\hdots,p_m\big\}$ of Eq.(\ref{logistica}), is to set a ``switching'' control, as follows 

\begin{equation}\label{DFCorbit}
u_k\!=\!\left\{\!\begin{array}{cc }
\!\gamma_i\,(x_k\!-\!x_{k-m}), &\ {\rm if }\quad  \sqrt{\sum\limits^m_{j=0}|x_{k-j}\!-\!p_{(i-j)\,{\rm mod}(m)}|^2}\ls\frac{\ts \varep}{\sqrt{2}},\ 1\!\ls\!i\!\ls\!m\, ,\\
 0  &  {\rm otherwise}\, ,
\end{array}\right.
\end{equation}
provided that $0\!<\!\varep\!<\!\frac{\ts \| P_i\!-\!P_j\|}{\ts \sqrt{2}},\ \forall\, i\!\neq\!j$, being $P_i\!=\!(p_i,p_{i-1},\hdots,p_1,p_m,p_{m-1},\hdots,p_i),\ 1\!\ls\! i\!\ls\! m.$
\medskip

The system (\ref{map-control})--(\ref{DFCorbit}) yields to a $(m\!+\!1)$-dimensional system (with $m$ switches), with $({\rm x}^1,\hdots,{\rm x}^{m+1})$ as its state variables:
\begin{equation}\label{syst_m+1}
\left\{\begin{array}{cl}
{\rm x}^1_{k+1}=&f({\rm x^1_k},r_0\!+\!\gamma_i\,({\rm x}^1_k\!-\!{\rm x}^{m+1}_k))\, ,\qquad {\rm if}\quad \|({\rm x}^1_k,\hdots,{\rm x}^{m+1}_k)-P_i\|\ls\frac{\ts \varep}{\sqrt{2}}\\
{\rm x}^2_{k+1}=&{\rm x}^1_k\, , \\
{\rm x}^3_{k+1}=&{\rm x}^2_k\, , \\
\vdots  & \vdots \\
{\rm x}^{m+1}_{k+1}=&{\rm x}^m_k\, .
\end{array}\right.
\end{equation}


Note that $\big\{P_i\!=\!(p_i,p_{i-1},\hdots,p_1,p_m,p_{m-1},\hdots,p_i)\big\}_{1\ls i\ls m}$ is $m$-UPO of the free system ($\gamma_i\!=\!0,\ \forall\, i$), and it is preserved when (\ref{DFCorbit}) is applied. 

The Jacobian matrix of the system (\ref{syst_m+1}) at $P_i$ is a $(m\!+\!1)\!\times\!(m\!+\!1)$--companion matrix given by
\begin{equation}
J_i=\left(\begin{array}{ccccc}
a^{(i)}_{11} & 0 & \cdots & 0 & a^{(i)}_{1(m+1)}\\
1 & 0 & \cdots & 0 & 0 \\
0 & 1 & \cdots &0 & 0 \\
\vdots & \vdots & \ddots & \vdots & \vdots \\
0 & 0 & \cdots & 1 & 0
\end{array}\right)
\end{equation}
where 
\begin{eqnarray*}
a^{(i)}_{11}&= &f_x(p_i,r_0)+\gamma_i\,f_r(p_i,r_0)\\ 
a^{(i)}_{1(m+1)}&=& -\gamma_i\,f_r(p_i,r_0)
\end{eqnarray*}
and, in particular for the logistic case, they become
\begin{eqnarray*}
a^{(i)}_{11}&= &\gamma_i\,p_i(1-p_i)+r_0(1-2p_i)\\ 
a^{(i)}_{1(m+1)}&=& -\gamma_i\,p_i(1-p_i)
\end{eqnarray*}
If each $\gamma_i$ could be chosen so that each $J_i$ has all its eigenvalues of modulus less than one, then, UPO stabilization arises as for SPF modulation in Proposition 1. Unfortunately, for $m\!>\!1$, a range of values for $\gamma_i$ to fulfill this requirement does not exist in most cases. Even for the $2$-UPO of the logistic map, there is no solution.

Indeed, the conditions for the orbit to be a.s. involve the product of the matrices, $\prod\limits_{i=1}^mJ_i$. Namely it should be found $\gamma_i$ making this product have all eigenvalues of modulus less than one. Note that this condition may be seen as the analogous of (\ref{prod_deriv}) for the control (\ref{control_mod}). Moreover, it reduces to the (only) Jacobian matrix of the system, for the fixed point stabilization (Eq. (\ref{DFCptofijo})).

As each $J_i$ is a companion matrix, the characteristic polynomial, $\chi(\lambda)$, of the product $\prod\limits_{i=1}^mJ_i$ is obtained by using tools from \cite{KV} and it results
\begin{equation*}
\chi(\lambda)=- \lambda^{m+1}+ \prod_{i=1}^m \big( a^{(i)}_{1(m+1)}+\lambda\, a_{11}^{(i)}\big)\, , 
\end{equation*}
so, for the logistic case,
\begin{equation}\label{pol-char}
 \chi(\lambda)= - \lambda^{m+1}+ \prod_{i=1}^m \big(-\gamma_i\,p_i(1-p_i)+\lambda(\gamma_i\,p_i(1-p_i)+r_0(1-2p_i))\big)\, .
\end{equation}
Therefore, the control is proposed as:
\begin{equation}\label{XXX}
u_k\!=\!\left\{\!\begin{array}{cc}
 0  & \ {\rm if}\quad k<k_0\, ,\\
\!\gamma_{(i+k-k_0)\,{\rm mod}(m)}\,(x_k\!-\!x_{k-m}), &\ {\rm if}\quad  k \gs k_0\, ,
\end{array}\right.
\end{equation}
provided that $0\!<\!\varep\!<\!\frac{\ts \| P_i\!-\!P_j\|}{\ts \sqrt{2}},\ \forall\, i\!\neq\!j$, being $k_0$=$k_0(x_0,\varep)$=$\min\{k\!\gs\!m:\|(x_{k+j},\hdots,x_{k+j-m})\!-\!P_{i+j}\|\ls \frac{\ts \varep}{\ts \sqrt{2}},\ {\rm for\ some}\ i,\ 1\!\ls\! i\!\ls\! m \}$. 
Existence of $k_0$ is a consequence of ergodicity. Introducing
$$X_k=\left(\begin{array}{l}{\rm x}^1_k \\ {\rm x}^2_k\\ \ \vdots\\ {\rm x}^{m+1}_k \end{array}\!\!\right) \qquad {\rm and} \qquad F_{\gamma_i}(X_k)=\left(\begin{array}{c} f({\rm x^1_k},r_0+\gamma_i({\rm x}_k^1-{\rm x}_k^{m+1})) \\ \vdots \\ {\rm x}_k^m  \end{array}\right)\ ,$$
the controlled system (\ref{map-control})--(\ref{XXX}) is equivalent to 
$$X_{k+1}=F_{\gamma_{(i+k-k_0)\,{\rm mod}(m)}}(X_k) \quad \forall\, k\!\gs\! k_0\ ,$$
which has $\{P_i\}_{1\ls i\ls m}$ as $m$-UPO. Note that
$$D_XF_{\gamma_i}(P_i)=J_i\, ,$$
and $k_0\!=\!\min\{k\!\gs\!m:\|X_k\!-\!P_i\|\!\ls\!\frac{\ts \varep}{\sqrt{2}}\}$. Next, the success of the proposed scheme is proved following the same steps that in the proof of Proposition 2:
\medskip

{\bf Proposition 3:} {\it Let the controlled system (\ref{map-control})--(\ref{XXX}) with $\gamma_j,\ 1\!\ls\! j\!\ls\! m$, such that all the roots of $\chi(\lambda)$ are of modulus less than one, and $\gamma=\max\limits_{1\ls i\ls m}|\gamma_j|$. Then, there exists $\varep_0,\ 0\!<\!\varep_0\!<\!\min\limits_{i\neq j}\Big\{\frac{\ts \|P_i\!-\!P_j\|}{\ts \sqrt{2}};\frac{\ts \delta}{\ts \gamma}\Big\}$ such that for all $\varep,\ 0\!<\!\varep\!<\!\varep_0$ and for almost every initial condition $x_0\!\in\!\cal{B}$, it verifies $|u_k|\!\ls\!\delta$, and that $(x_k)_{k\gs 1}$ converges to the $m$-UPO $\{p_1,\hdots,p_m\}$. }

\ul{Proof:} As the roots of $\chi(\lambda)$ are of modulus less than one, then $\big\|\prod\limits_{i=1}^mJ_i\big\|\!<\!1$. Let $\varep'_0\!=\!\min\limits_{i\neq j}\Big\{\frac{\ts \|P_i\!-\!P_j\|}{\ts \sqrt{2}}\Big\}$ and $0\!<\!\sigma\!<\!1$ such that
$$\big\|\prod\limits_{i=1}^mJ_i\big\|<\sigma <1$$
which is the same as
$$\big\|D_XF_{\gamma_1}(P_1)D_XF_{\gamma_2}\big(F_{\gamma_1}(P_1)\big)\cdots D_XF_{\gamma_m}\big(F_{\gamma_{m-1}}\big(\cdots\big(F_{\gamma_1}(P_1)\big)\big)\big) \big\|<\sigma\ .$$
Then, there exists $\varep_1,\ 0\!<\!\varep_1\!<\!\varep'_0$ such that, for $X: \|X\!-\!P_1\|\!\ls\!\frac{\ts \varep_1}{\ts \sqrt{2}}$
$$\big\|D_XF_{\gamma_1}(X)D_XF_{\gamma_2}\big(F_{\gamma_1}(X)\big)\cdots D_XF_{\gamma_m}\big(F_{\gamma_{m-1}}\big(\cdots\big(F_{\gamma_1}(X)\big)\big)\big) \big\|\ls\sigma\ .$$
We can proceed analogously for $P_2,\hdots,P_m$, and obtain $\varep_2,\hdots,\varep_m$. Let $\varep_0\!=\!\min\limits_{1\ls i\ls m}\varep_i$, and fix $\varep\!<\!\varep_0$. Without loss of generality, we assume that condition $\|X_{k_0}\!-\!P_1\|\!\ls\!\frac{\ts \varep}{\ts \sqrt{2}}$ is verified, for the first time, by $i\!=\!1$. Then
\begin{eqnarray*}
\|X_{k_0+m}-P_1\|&=&\big\| F_{\gamma_m}\circ F_{\gamma_{m-1}}\circ\cdots\circ F_{\gamma_1}(X_{k_0})- F_{\gamma_m}\circ F_{\gamma_{m-1}}\circ\cdots\circ F_{\gamma_1}(P_1)\big\|\\
     &\ls& \big\| D_X\big(F_{\gamma_m}\circ F_{\gamma_{m-1}}\circ\cdots\circ F_{\gamma_1}\big)(\xi_1)\big\| \, \|X_{k_0}-P_1\| \\
               &\ls& \sigma \, \frac{\ts \varep}{\ts \sqrt{2}} < \frac{\ts \varep}{\ts \sqrt{2}}\ .
\end{eqnarray*}
where $\xi_1$ is a point within the segment joining $P_1$ and $X_{k_0}$. By the same arguments it is obtained that
\begin{equation}\label{convergence}
\|X_{k_0+j+\ell m}-P_{j+1}\|\ls\sigma^{\ell} \, \frac{\ts \varep}{\ts \sqrt{2}} < \frac{\ts \varep}{\ts \sqrt{2}}\qquad j\!=\!1,\hdots,m\!-\!1,\ \ell\!\gs\!1\ .
\end{equation}   
Hence, the convergence of $(X_k)_{k\gs 1}$ to the $m$-UPO $\{P_i\}_{1\ls i\ls m}$ follows, and as a consequence, the convergence of $(x_k)_{k\gs 1}$ to the $m$-UPO $\{p_1,\hdots,p_m\}$. 

From Eq. (\ref{convergence}) it also results $$|x_{k_0+j+\ell m}-x_{k_0+j+m(\ell-1)}|<\varep\qquad j\!=\!1,\hdots,m\!-\!1,\ \ell\!\gs\!1\ ,$$
and so, the control bounds are obtained. \hfill{$\square$}
\medskip

It remains to solve the problem of finding a range of $\gamma_i$'s for orbit stability, i.e., for the roots of $\chi(\lambda)$ (Eq. (\ref{pol-char})) to be within the unit circle. In the general situation, it is not possible to obtain the $\gamma_i$'s explicitely, neither to ensure a range finding. However, we first attempt to look for a solution, by fixing one $\gamma_i\!\neq\!0$ and setting the rest $\gamma_j\!=\!0,\ \forall\,j\!\neq\!i$, in order to simplify the equation. In this way, we obtain $\lambda\!=\!0$ as a $(m\!-\!1)$-multiplicity root of $\chi(\lambda)$, and the other two being the roots of a quadratic function, that is, 
\begin{equation}\label{poli}
\chi(\lambda)=\lambda^{m-1}\big(\lambda^2+A\lambda +B_i\big)\, ,
\end{equation}
where, $A=\prod\limits_{j=1}^m a_{11}^{(j)}$, $B_i=a^{(i)}_{1(m+1)}C_i$, and $C_i=\prod\limits_{j\neq i}a_{11}^{(j)}$. A range $[\gamma_i^{\inf},\gamma_i^{\sup}]$ should be obtained, from the conditions on $A$ and $B_i$ that make the roots of the quadratic function be within the unit circle (\cite{K}). 

For the logistic map, it results 
\renewcommand{\arraystretch}{1.5}
\begin{equation}\label{gama-inf-sup}
\begin{array}{lll}
\gamma_i^{\inf}=\left\{\begin{array}{cl}\frac{-1\!-\!r\,(1\!-\!2p_i)C_i}{2p_i(1\!-\!p_i)C_i} & {\rm if}\ C_i>0 \\
\frac{1}{p_i(1\!-\!p_i)C_i}& {\rm if}\ C_i<0\end{array}\right. & {\rm and} & 
\gamma_i^{\sup}=\left\{\begin{array}{cl}\frac{1}{p_i(1\!-\!p_i)C_i}& {\rm if}\ C_i>0 \\
\frac{-1\!-\!r\,(1\!-\!2p_i)C_i}{2p_i(1\!-\!p_i)C_i}& {\rm if}\ C_i<0\end{array}\right.
\end{array}
\end{equation} \renewcommand{\arraystretch}{1.2}
whenever $\gamma_i^{\inf}<\gamma_i^{\sup}$, which implies the following condition on $r, p_1,\hdots,p_m$ to be fullfiled:
\begin{equation}\label{condition}
\left|1+r^m\prod_{j=1}^{m}\,(1\!-\!2p_j)\right|<2\,.
\end{equation}
This procedure is applied to the $4$-UPO $\{p_1,p_2,p_3,p_4\}$ of Eq. (\ref{logistica}). 
For the parameter value $r\!=\!3.62$, $p_1\!\approx\!0.5522$, $p_2\!\approx\!0.8951$, $p_3\!\approx\!0.3398$, $p_4\!\approx\!0.8121$, the condition (\ref{condition}) is accomplished. Taking, for example, $\gamma_1\!=\!\gamma_2\!=\!\gamma_3\!=\!0$,  a range for $\gamma_4$ is obtained from (\ref{gama-inf-sup}). Figure 4 shows this UPO being stabilized by the control (\ref{XXX}), for initial condition $x_0\!=\!0.5$ and $\varep\!=\!0.05$. The performance of the control $u_k$ is shown in the same figure.

Unfortunately, this methodology does not work for higher order period orbits anymore, and, even for $m\!=\!4$, it cannot be applied for parameter values above $r\!\approx\!3.625$ because the Eq. (\ref{condition}) is not accomplished. However, for values of $r$ over $3.625$, by means of a detailed heuristic search on the coefficient values of $\chi(\lambda)$, 4-tuples $[\gamma_1,\gamma_2,\gamma_3,\gamma_4]$ can be found such that all the roots of $\chi(\lambda)$ are within the unit circle, and so, stabilizing the $4$-UPO (see Fig. 5). This shows that the control strategy (\ref{XXX}) overcomes the Pyragas method, by  which the 4-cycle logistic map control can not be maintained above $r\!\approx\!3.62$ (see \cite{Socolar} p.50). Namely, by control (\ref{XXX}), the 4-UPO is stabilized for $3.625\!\ls\!r\!\ls\!3.67$ and, even in presence of noise (Fig. 6). 

\subsection{Improving EDFC modulation}

The extension of the Pyragas method, called ``extended time delay autosynchronization system'' (ETDAS) in \cite{Socolar}, and ``extended delayed feedback control'' (EDFC) in \cite{Ko}, is defined as
\begin{equation}\label{etdas}
u_k=\gamma\,(x_k\!-\!x_{k-m})+R\,u_{k-m}\, ,
\end{equation}
with $0\!\ls\!R\!<\!1$, so the case $R\!=\!0$ reduces to TDAS. It is introduced to stabilize higher instabilities. Indeed, it is reported in \cite{Socolar} that with $R\!=\!0.5$, stabilization of the 4-UPO can be maintained up to the parameter value $r\!\approx\!3.75$. With these results in mind, we have designed a kind of extension like (\ref{etdas}), but based on the  control (\ref{XXX}), that consists of replacing it by
\begin{equation}\label{XXXX}
u_k\!=\!\left\{\!\begin{array}{cc}
 0  & \ {\rm if}\quad k<k_0\, ,\\
\!\gamma_{(i+k-k_0)\,{\rm mod}(m)}\,(x_k\!-\!x_{k-m})+R\,u_{k-m} , &\ {\rm if}\quad  k \gs k_0\, ,
\end{array}\right.
\end{equation}
Now, the system (\ref{map-control})--(\ref{XXXX}) yields to a $2m$-dimensional system (with $m$ switches), with $({\rm x}^1,\hdots,{\rm x}^m,{\rm u}^1,\hdots,{\rm u}^m)$=$(x_k,\hdots,x_{k-m+1},u_k,\hdots,u_{k-m+1})$ as its state variables. For the logistic map, the $m\!=\!4$ case is given by
\renewcommand{\arraystretch}{1}
\begin{equation}\label{syst_2m}
\left\{\begin{array}{cll}
{\rm x}^1_{k+1}&=&(r\!+\!{\rm u^1_k})\,{\rm x}^1_k\,(1\!-\!{\rm x}^1_k)\, ,\qquad {\rm if}\quad \|({\rm x}^1_k,\hdots,{\rm x}^m_k)-{\widetilde P}_i\|\ls\frac{\ts \varep}{\sqrt{2}}\\
{\rm x}^2_{k+1}&=&{\rm x}^1_k\, , \\
{\rm x}^3_{k+1}&=&{\rm x}^2_k\, , \\
{\rm x}^4_{k+1}&=&{\rm x}^3_k\, , \\
{\rm u}^1_{k+1}&=& \gamma_i\,[(r\!+\!{\rm u^1_k})\,{\rm x}^1_k\,(1\!-\!{\rm x}^1_k)-{\rm x}^4_k]+R\,{\rm u^4_k}\, ,\\
{\rm u}^2_{k+1}&=&{\rm u}^1_k\, , \\
{\rm u}^3_{k+1}&=&{\rm u}^2_k\, , \\
{\rm u}^4_{k+1}&=&{\rm u}^3_k\, .
\end{array}\right.
\end{equation}
Note that $\big\{{\widetilde P}_1\!=\!(p_1,p_4,p_3,p_2,0,0,0,0)$, ${\widetilde P}_2\!=\!(p_2,p_1,p_4,p_3,0,0,0,0)$, ${\widetilde P}_3\!=\!(p_3,p_2,p_1,p_4,0,0,0,0)$, ${\widetilde P}_4\!=\!(p_4,p_3,p_2,p_1,0,0,0,0)\big\}$ is a $4$-UPO of the free system ($\gamma_i\!=\!0,\ \forall\, i$), and it is preserved when (\ref{XXXX}) is applied. 

The Jacobian matrix of the system (\ref{syst_2m}) at ${\widetilde P}_i$ is, in this case, a $8\!\times\!8$ matrix given by
\renewcommand{\arraystretch}{0.8}
\begin{equation} 
J_i=\left(\begin{array}{cccrcccc}
a(i) & 0 & 0 & 0 & b(i)& 0 & 0 & 0 \\
 1 & 0 & 0 & 0 & 0 & 0 & 0 & 0 \\
 0 & 1 & 0 & 0 & 0 & 0 & 0 & 0 \\
 0 & 0 & 1 & 0 & 0 & 0 & 0 & 0 \\
c(i) & 0 & 0 & -\gamma_i & d(i) & 0 & 0 & R \\
 0 & 0 & 0 & 0 & 1 & 0 & 0 & 0 \\
 0 & 0 & 0 & 0 & 0 & 1 & 0 & 0 \\
 0 & 0 & 0 & 0 & 0 & 0 & 1 & 0 
 \end{array}\right)
 \end{equation}
where $a(i)\!=\!r\,(1\!-\!2\,p_i)$, $b(i)\!=\!p_i\,(1\!-\!p_i)$, $c(i)\!=\!\gamma_i\,r\,(1\!-\!2\,p_i)$ and $d(i)\!=\!\gamma_i\,p_i\,(1\!-\!p_i)$. As $c(i)\!=\!\gamma_i\,a(i)$ and $d(i)\!=\!\gamma_i\,b(i)$, the characteristic polynomial is given by
\begin{equation}\label{poli}
\chi(\lambda)=\lambda^3\,\big[\lambda^5-(a(i)+\gamma_i\,b(i))\lambda^4-R\,\lambda+a(i)R+\gamma_i\,b(i)\big]
\end{equation}
Applying the heuristic algorithm to search on the coefficient values of $\chi(\lambda)$, we find four-tuples, $[\gamma_1,\gamma_2,\gamma_3,\gamma_4]$ and can stabilize the 4-UPO for $3.67\!<\!r\!\ls\!3.8$. Thus, the range of parameter values for which the 4-cycle-logistic-map is stabilized, is enlarged. The Figures \ref{7} and \ref{8} illustrate the orbit stabilization and the performance of the corresponding control for some values of $r$.

\section{Final Remarks}

We just state  a few comments about these developments when applied to $n$-dimensional systems.

For SPF modulation, Eq. (\ref{controlbetas}) remains all the same save that control gain is $\beta_i^\top\!\in\!\Bbb{R}^n$ and the distance is given by the Euclidean norm in $\Bbb{R}^n$. The Jacobian of the controlled system in each $p_i$, $1\!\ls\! i\!\ls\!m$, is given by:
$$\frac{d\,f_{\beta_i}}{dx}(p_i)=\frac{\part f}{\part x}(p_i,r_0)+\frac{\part f}{\part r}(p_i,r_0)\,\beta_i^\top$$
where $A_x\!=\!\frac{\part f}{\part x}(p_i,r_0)$ and $A_r\!=\!\frac{\part f}{\part r}(p_i,r_0)$ are $n\!\times\!n$ and $n\!\times\!1$ matrices, respectively. Then, the condition for stability requires that $A_x+A_r\beta_i^\top$ must have the whole of its eigenvalues of modulus less than one. A sufficient condition for the existence of each $\beta_i^\top\!\in\!\Bbb{R}^n$ is that rank $\{A_r,A_xA_r,\hdots,A_x^{n-1}A_r\}\!=\!n$ and pole-placement methods is a key tool to obtain $\beta_i$ \cite{Ogata}. In turn, this implies that $\left\|\frac{df_{\beta_i}}{dx}(p_i)\right\|\!<\!\sigma$ (with the matrix norm induced by the Euclidean one).
Defining $\beta\!=\!\max\limits_{1\ls i\ls m}\|\beta_i\|$, the results of Section 2 are valid in the $n$-dimensional case.

Analogously, for the DFC modulation, a vector $\gamma_i^\top\!\in\!\Bbb{R}^n$ becomes the control gain for the $n$-dimensional version of Eq. (\ref{DFCorbit}). The stability criterion involves the product $\prod\limits_{i=1}^m J_i$, being $J_i$ a $n(m\!+\!1)\!\times\!n(m\!+\!1)$ matrix given by:
\begin{equation}
J_i=\left(\begin{array}{ccccc}
A^{(i)}_{11} & 0 & \cdots & 0 & A^{(i)}_{1(m\!+\!1)}\\
I_n & 0 & \cdots & 0 & 0 \\
0 & I_n & \cdots &0 & 0 \\
\vdots & \vdots & \ddots & \vdots & \vdots \\
0 & 0 & \cdots & I_n & 0
\end{array}\right)
\end{equation}
where $A_{11}^{(i)}=\frac{\part f}{\part x}(p_i,r_0)+\frac{\part f}{\part r}(p_i,r_0)\, \gamma_i^\top$ and $A_{1\,(m+1)}^{(i)}=-\frac{\part f}{\part r}(p_i)\, \gamma_i^\top$ are both $n\!\times\!n$ matrices and $I_n$ is the $(n\!\times\!n)$-identity matrix.

In the same way, the EDFC algorithm of Eq. (\ref{etdas}) may also be formulated for the $n$-dimensional case, its stability analysis yielding to a $2mn$-dimensional system (with $m$ switches).

As in the SPF modulation, controllability tools become useful to obtain adequate control gains. Defining $\gamma\!=\!\max\limits_{1\ls i\ls m}\|\gamma_i\|$ the arguments to prove the validity of these methods are straightfull generalized from Section 3; although full description of them become tiresome. In spite of this, the implementation (and the corresponding analysis) on specific systems like [Henon]...


\newpage

\begin{figure}[h]
\hspace{-3cm}
\includegraphics[width=21cm,height=11cm]{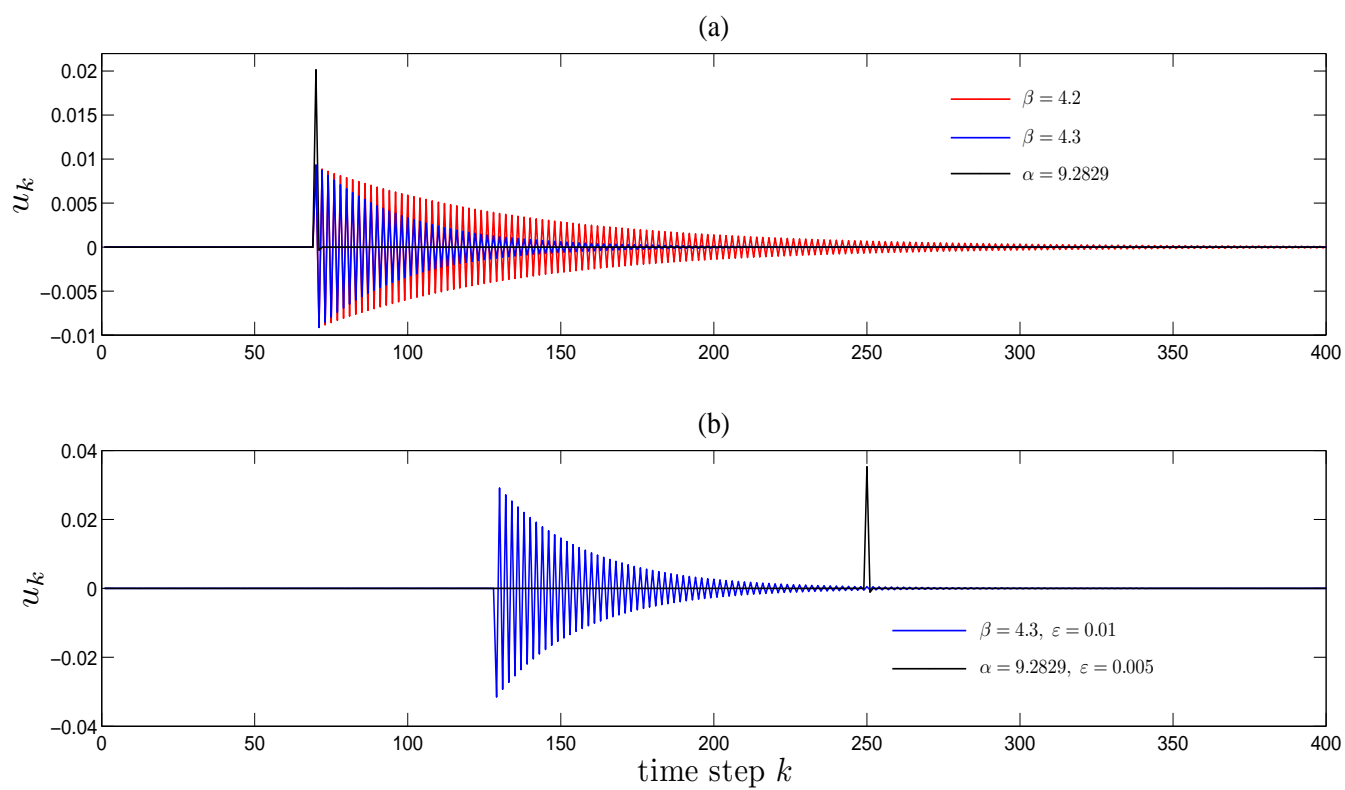}
\caption{\small (a) Performance of the controls (\ref{controlbetas}) and (\ref{DFCptofijo}) to stabilize the Eq. (\ref{logistica}) with $r_0\!=\!3.8$, about the fixed point $p\!\approx\!0.736$, $x_0\!=\!0.94,\varep\!=\!0.005$, ($\beta^{\inf}\!\approx\!4.126,\beta^{\sup}\!\approx\!14.44$). (b) Idem (a), but keeping the control effort, and varying the $\varep$ values, $x_0\!=\!0.5$.}
\label{2}
\end{figure}

\begin{figure}[h]
\hspace{-3cm}
\includegraphics[width=21cm,height=11cm]{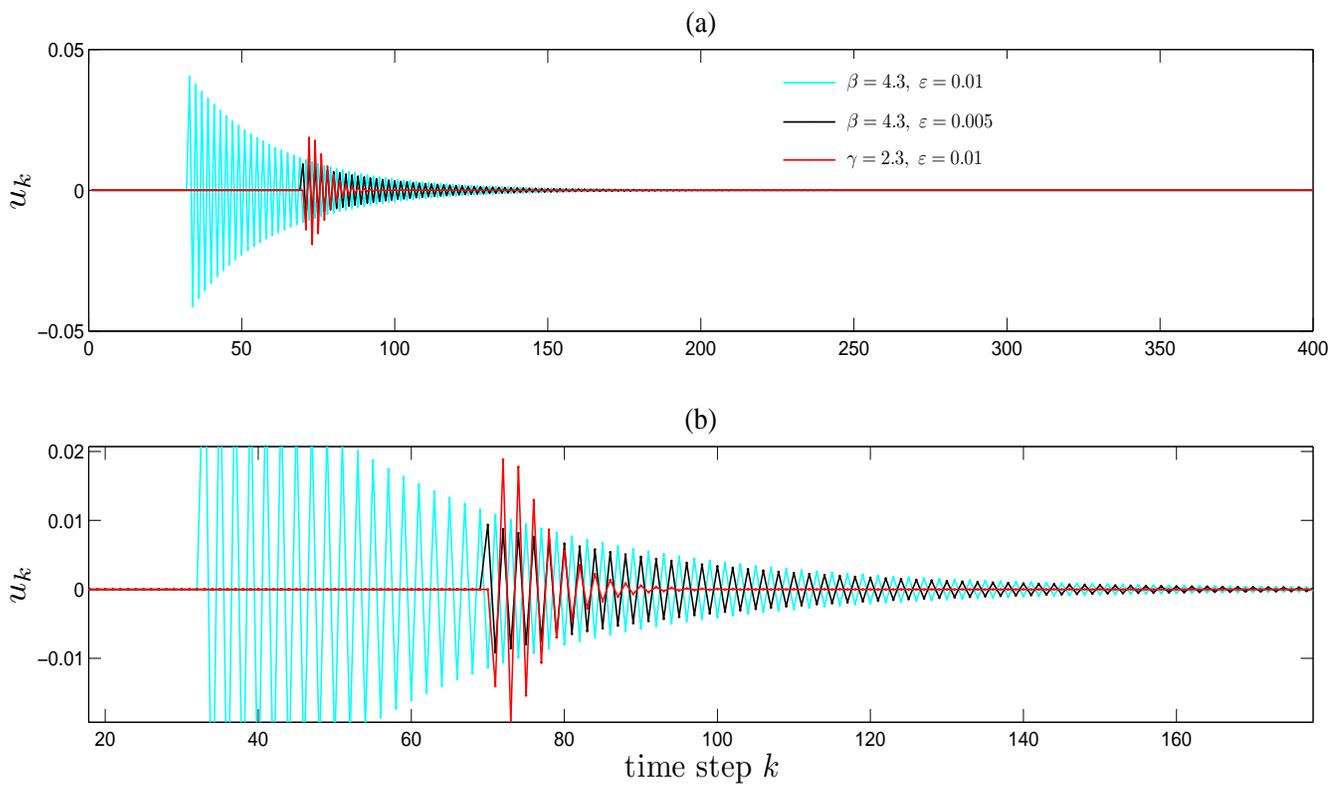}
\caption{\small (a) Performance of the controls (\ref{controlbetas}) and (\ref{DFCptofijo}) to stabilize the Eq. (\ref{logistica}) with $r_0\!=\!3.8$, about the fixed point $p\!\approx\!0.736$, $x_0\!=\!0.94,\ y_0\!=\!0.45$. (b) Zoom in of (a).}
\label{3}
\end{figure}


\begin{figure}[h]
\hspace{-3cm}
\includegraphics[width=21cm,height=11cm]{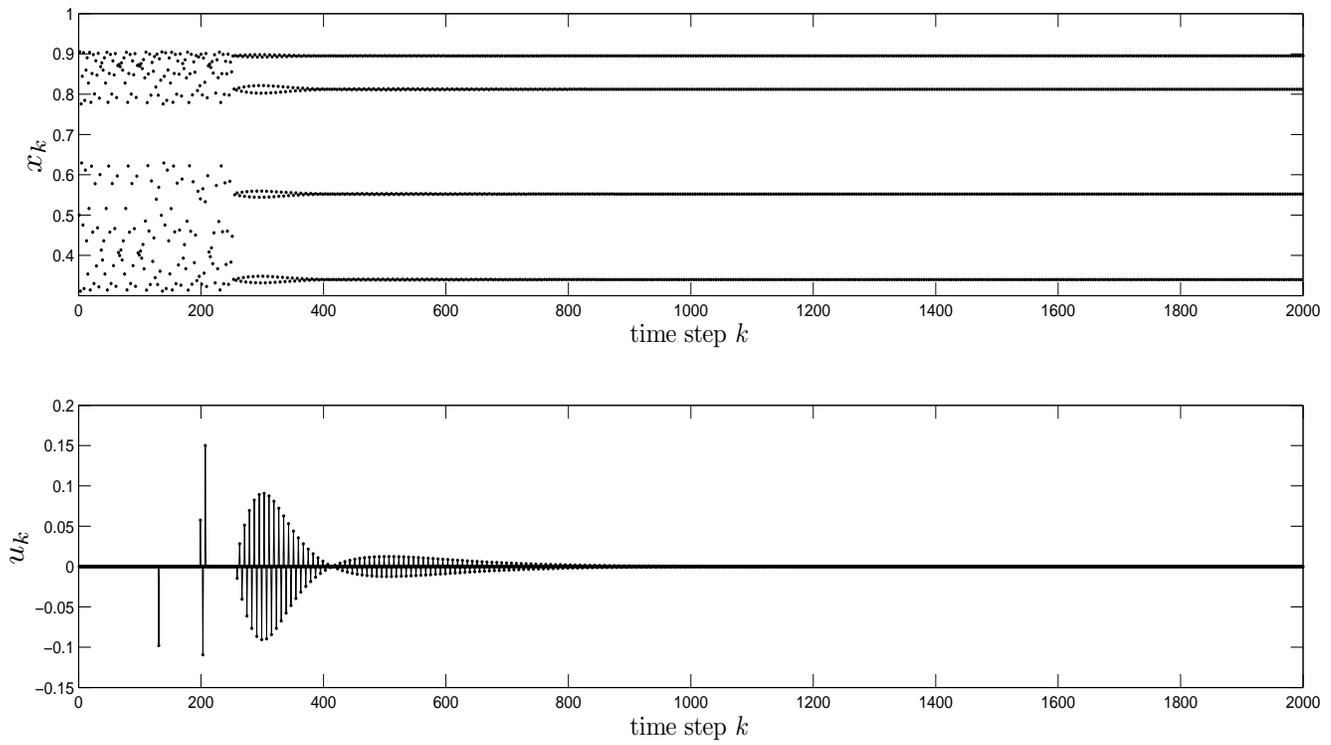}
\caption{\small States $x_k$ and controls $u_k$ of Eq. (\ref{XXX}) with $\gamma_1\!=\!\gamma_2\!=\!\gamma_3\!=\!0$, $\gamma_4\!=\!4.7997$, $\varep\!=\!0.05$, applied to stabilize $4$-UPO of Eq. (\ref{logistica}) for $r\!=\!3.62$.}
\label{4}
\end{figure}

\begin{figure}[h]
\hspace{-3cm}
\includegraphics[width=21cm,height=11cm]{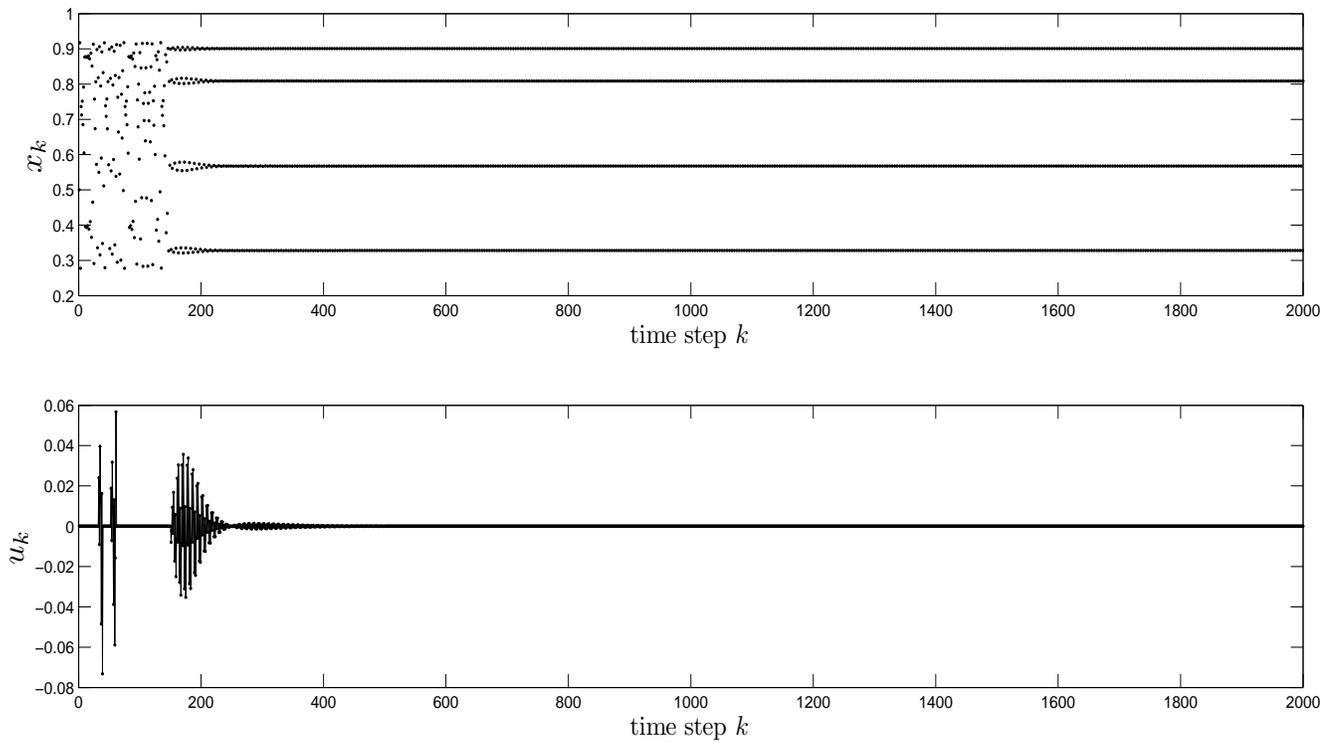}
\caption{\small States $x_k$ and controls $u_k$ of Eq. (\ref{XXX}), with $\gamma_1\!=\!0.4,\gamma_2\!=\!4.97156,\gamma_3\!=\!-0.598,\gamma_4\!=\!2.09$, $\varep\!=\!0.05$, applied to stabilize $4$-UPO of Eq. (\ref{logistica}) for $r\!=\!3.67$.}
\label{5}
\end{figure}

\begin{figure}[h]
\hspace{-3cm}
\includegraphics[width=21cm,height=11cm]{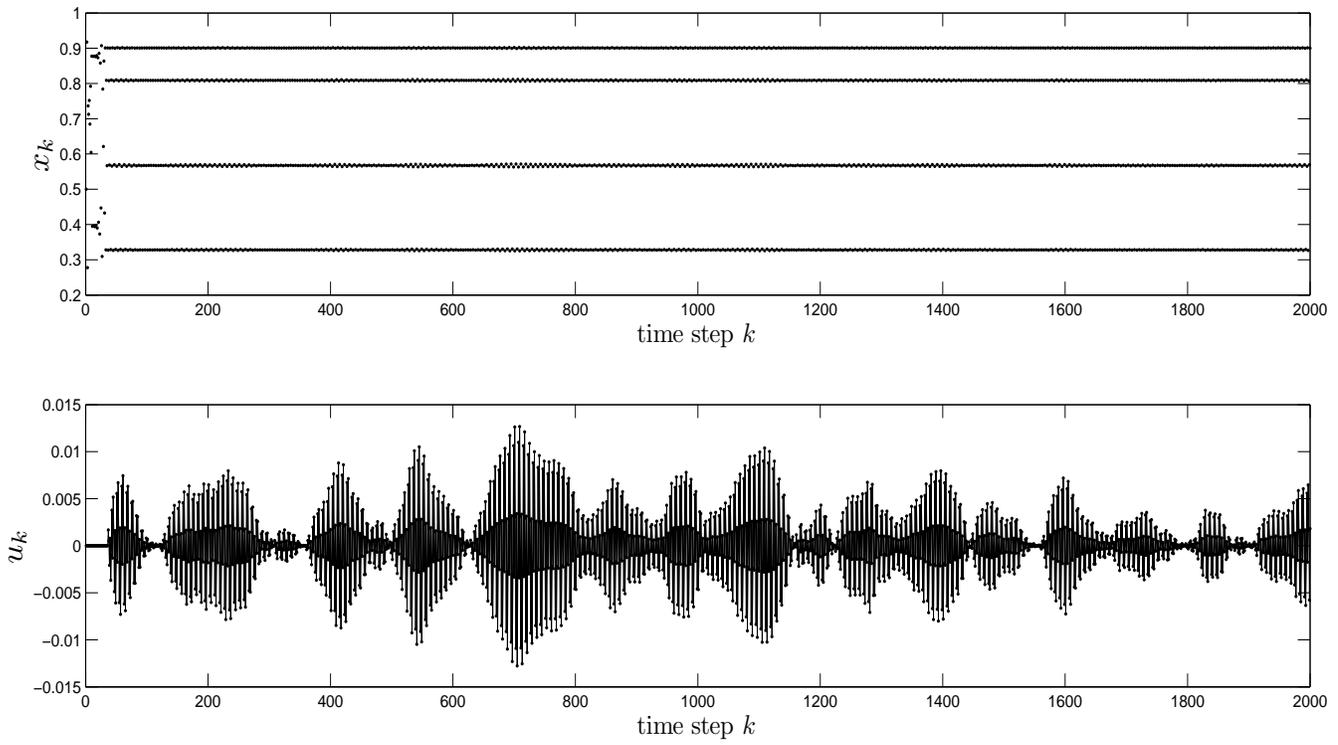}
\caption{\small States $x_k$ and controls $u_k$ of Eq. (\ref{XXX}), with $\gamma_1\!=\!0.4,\gamma_2\!=\!4.97156,\gamma_3\!=\!-0.598,\gamma_4\!=\!2.09$, $\varep\!=\!0.05$, applied to stabilize $4$-UPO of Eq. (\ref{logistica}) for $r\!=\!3.67$, under the effect of additive noise modeled by $2\!\times\!10^{-5}\sigma_k$, $\sigma_k\!\sim\!N(0,1)$.}
\label{6}
\end{figure}

\begin{figure}[h]
\hspace{-3cm}
\includegraphics[width=21cm,height=11cm]{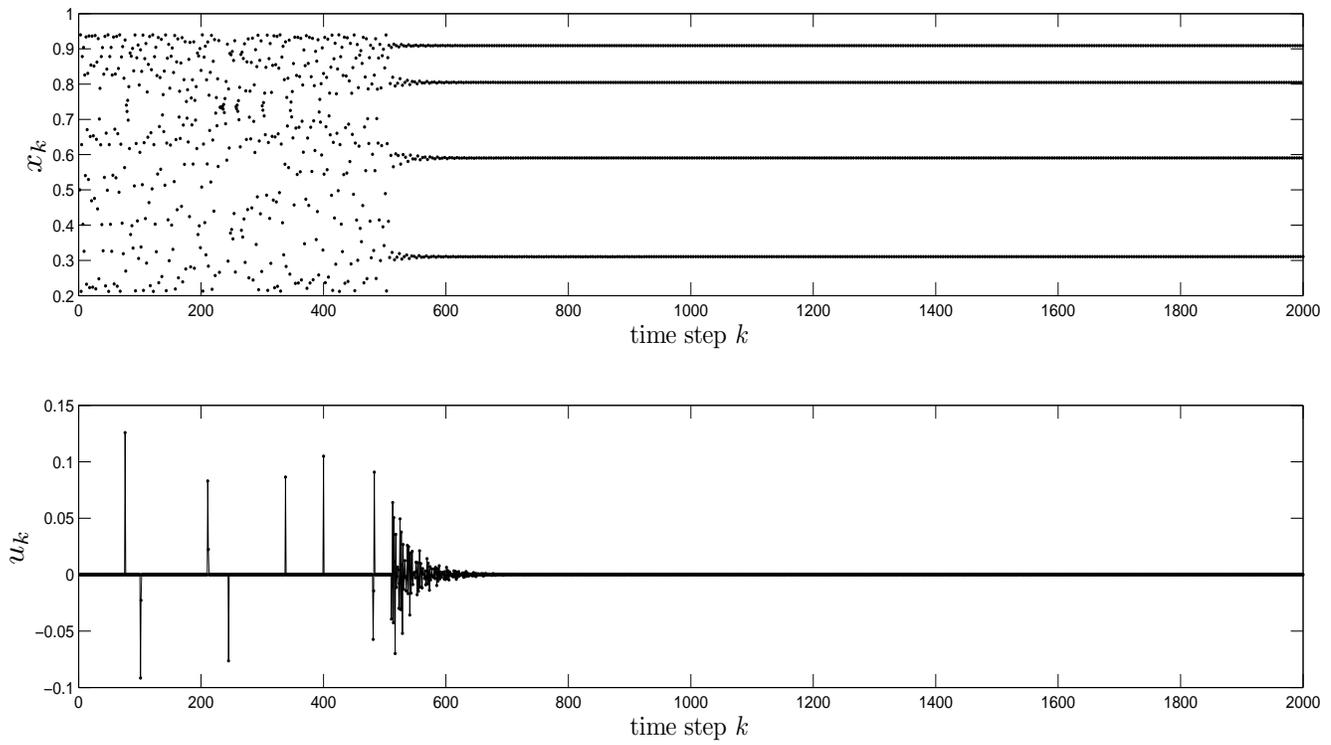}
\caption{\small States $x_k$ and controls $u_k$ of Eq. (\ref{XXXX}), with $\gamma_1\!=\!1.333,\gamma_2\!=\!6.79,\gamma_3\!=\!-0.6999,\gamma_4\!=\!3.601$, $R\!=\!0.3$, $\varep\!=\!0.05$, applied to stabilize $4$-UPO of Eq. (\ref{logistica}) for $r\!=\!3.76$.}
\label{7}
\end{figure}

\begin{figure}[h]
\hspace{-3cm}
\includegraphics[width=21cm,height=13cm]{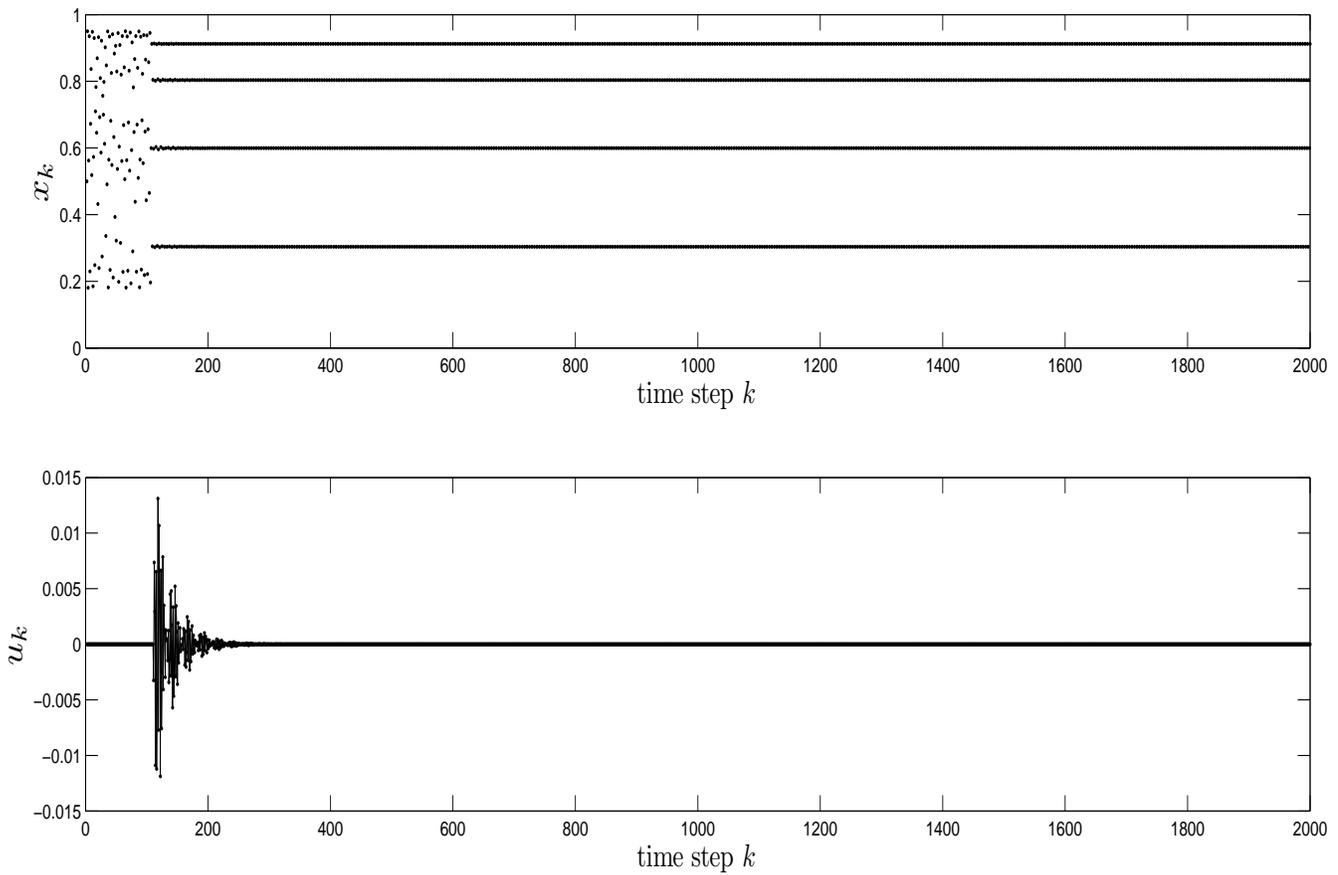}
\caption{\small States $x_k$ and controls $u_k$ of Eq. (\ref{XXXX}), with $\gamma_1\!=\!3.50293,\gamma_2\!=\!1.38,\gamma_3\!=\!7.49498,\gamma_4\!=\!-1.181$, $R\!=\!0.3$, $\varep\!=\!0.05$, applied to stabilize $4$-UPO of Eq. (\ref{logistica}) for $r\!=\!3.8$.}
\label{8}
\end{figure}

\end{document}